\documentclass[11pt]{article}
\usepackage{latexsym}
\title{Three-manifolds of positive Ricci curvature and convex weakly umbilic boundary}
 
\author{Jean C. Cortissoz}
\date{}

\newtheorem{thm}{Theorem}[section]
\newtheorem{prop}[thm]{Proposition}
\newtheorem{lem}[thm]{Lemma}
\newtheorem{cor}[thm]{Corollary}

\newtheorem{theorem}{Theorem}
\newtheorem{conjecture}{Conjecture}

\newtheorem{defn}[thm]{Definition}

\newtheorem{rem}[thm]{Remark}
 
\begin{document}
\maketitle 
\begin{abstract}
 In this paper we study a boundary value problem in manifolds with
weakly umbilic boundary (the Second
Fundamental form of the boundary is a constant multiple of the metric). We show that if we start with
a metric of positive Ricci curvature and convex boundary (positive Second
Fundamental form), the flow uniformizes the curvature. 
 
\end{abstract}
 
\section{Introduction}\label{intro}

In this paper we consider the problem of deforming metrics in manifolds with boundary via
the Ricci flow. The Ricci flow in manifolds with boundary was first considered by Shen in
\cite{Shen}. The problem he considered was the following
\begin{equation}
\label{BVP}
\left\{
\begin{array}{l}
\partial_t g = -2Ric(g) \quad \mbox{in} \quad M\times\left[0,T\right)\\
h_g = \kappa g \quad \mbox{on} \quad\partial M \times\left[0,T\right)\\
g|_{t=0}=\overline{g}
\end{array}
\right.
\end{equation}
where $h_g$ is the Second Fundamental Form of the boundary 
with respect to the outward unit normal, and {\bf$\kappa$ is a constant}.
We also ask for the following compatibility condition (which is required in
order to show that the solution has enough regularity),
\[
h\left(\cdot,0\right)=\kappa \overline{g}\left(\cdot\right) .
\]

The proof of the basic short time existence result is 
given in \cite{Shen}. Here we state it
with a little more specifications which are necessary to justify certain applications of the 
Maximum Principle.

\begin{thm}
The Boundary Value Problem
(\ref{BVP})
has a unique continuous solution for a short time. This
solution is smooth in $\overline{M}\times \left(0,T\right)$,
and $Ric\left(g\right)$
is continuous in $\overline{M}\times [0,T)$.
\end{thm}

Also, in the same paper, the following result is proved

\begin{thm}[\cite{Shen}]
Let $\left(M,g\right)$ be a compact three-dimensional Riemannian manifold with
totally geodesic boundary and with positive Ricci curvature. Then $\left(M,g\right)$
can be deformed to $\left(M,g_{\infty}\right)$ via the Ricci flow such that
$\left(M,g_{\infty}\right)$ has constant positive curvature and totally geodesic boundary.
\end{thm}

Here we consider the normalized version of (\ref{BVP}). The normalization goes as follows.
Let $\psi\left(t\right)$ be such that for $\tilde{g}=\psi g$ we have $Vol\left(M\right)=1$. Then
we change the time scale by letting
\[
\tilde{t}\left(t\right)=\int_0^t \psi\left(t\right)\,d\tau ,
\]
then $\tilde{g}\left(\tilde{t}\right)$ satisfies the equation
\begin{equation}
\left\{
\label{BVP2}
\begin{array}{l}
\frac{\partial}{\partial \tilde{t}}\tilde{g}=\frac{2}{n}\tilde{r}\tilde{g}-2Ric\left(\tilde{g}\right) 
\quad \mbox{in} \quad M\times\left[0,T\right)\\
h_{g}=\kappa\left(\tilde{t}\right)\tilde{g} \quad\mbox{on} \quad \partial M \times\left[0,T\right)
\end{array}
\right.
\end{equation}
where 
\[
\tilde{r}=\int_{M} \tilde{R}\,dV .
\]
and $\tilde{R}$ is the scalar curvature of the metric $\tilde{r}$.

Our main result is
 
\begin{theorem}
\label{MainTheorem2}
If $\kappa\geq 0$ the solution 
of (\ref{BVP2}) exists for all time and converges exponentially to a metric
of constant sectional curvature and totally geodesic boundary.
\end{theorem}

We conjecture that
\begin{conjecture}
If $h=\kappa g$, for $\kappa$ any nonnegative function, then the metric $g_0$ can be 
deformed via the Ricci flow to a manifold of constant curvature and totally geodesic boundary.
\end{conjecture}

This paper is organized as follows. In Section \ref{MaximumPrinciples} we present
the basic Maximum Principle for tensors used in this work. In section \ref{Pinching},
we prove the basic pinching estimates, and we give an argument to prove Theorem \ref{MainTheorem2}.
Finally, in Section \ref{Derivative} we sketch a method to produce bounds on derivative
of the curvature up to the boundary from bounds in the curvature.

\section{Maximum Principle}
\label{MaximumPrinciples}

The Hopf Maximum Principle for the Ricci Flow is due to Shen (\cite{Shen}). We restate it
here in a slightly different way, very convenient for our purposes. First a definition.
\begin{defn}
Let $M_{ij}$ be a tensor and let $N_{ij}=p\left(M_{ij},g_{ij}\right)$ be
a polynomial in $M_{ij}$ formed by contracting products of $M_{ij}$ with itself
using the metric $g_{ij}$. We say that $N$ satisfies the \bfseries null eigenvector condition \normalfont
if whenever $v^i$ is a null eigenvector of $M_{ij}$, the we have $N_{ij}v^iv^j\geq0$ 

We say that $M$ satisfies the \bfseries normal derivative condition \normalfont 
at a point $p\in \partial M$, if for any null-eigenvector $v$ of $M_{ij}$,
\[
\left(M_{ij;\nu}\right)v^iv^j\geq 0
\]
\end{defn}
\begin{thm}[\cite{Shen}, \cite{Ham82}]
Let $(M,g)$ be a Riemannian manifold such that $R_{ij}\geq -\left(n-1\right)\omega g_{ij}$.
Suppose we have
\[
\frac{\partial}{\partial t}M_{ij}=\Delta M_{ij}+u^k\nabla_{k}M_{ij}+N_{ij}
\]
for some constant $\omega$, where $N=P\left(M_{ij},g_{ij}\right)$ satisfies the null-eigenvector
condition, and $M$ satisfies
the normal derivative condition. Then the condition $M_{ij}>0$ is preserved
under the flow.
\end{thm}

\section{Shen's Parabolic Simon Identity}

In this section we prove an important parabolic identity for the Second Fundamental form
of the boundary of $\left(M,g\right)$. This identity was derived by Shen 
(with what it seems to us some mistakes), in \cite{Shen} using Cartan's formalism 
for his computations.
We recast and rederive this formula using classical tensor notation. First we fix some notation.

\bfseries \textit{Notation. } \normalfont
The metric $g$ restricted to $\partial M$ will be denoted by $\overline{g}$. $\nabla$
will denote the connection of $g$, whereas $\overline{\nabla}$ will denote 
the (Levi-Civita) connection of $\overline{g}$. Covariant differentiation
with respect to $\nabla$ will be denoted by a vertical bar ($|$). Covariant
differntiation with respect to $\overline{\nabla}$ will be denoted by a semicolon ($;$).
 By $h_{\alpha\beta}$ we will
denote the Second Fundamental form of $\partial M$. Finally, greek numerals 
(except $\nu$, which we have chosen to represent the outward unit normal) denote quantities
in the boundary.

\hfill $\Box$

We are ready to establish,

\begin{prop}[Shen's Simon parabolic identity]
\begin{equation}
\begin{array}{rcl}
\frac{\partial}{\partial t}h_{\alpha\beta}
&=&\Delta h_{\alpha\beta}-H_{\alpha\beta}-R_{n\alpha\beta n|n}
  -HR_{n\alpha\beta n}\\
&&-g^{\gamma\delta}g^{\rho\omega}h_{\beta\rho}R_{\alpha\delta\gamma\omega}
  -g^{\gamma\delta}g^{\rho\omega}h_{\delta\rho}R_{\beta\gamma\alpha\omega}\\
&&-Hg^{\rho\omega}h_{\beta\rho}h_{\omega\alpha}+|A|^2h_{\alpha\beta}\\
&&-g^{\gamma\delta}g^{\rho\omega}\left(R_{\beta\delta\gamma\rho}h_{\omega\alpha}
  +R_{\beta\delta\alpha\rho}h_{\omega\gamma}\right)\\
&&-g^{\gamma\delta}g^{\rho\omega}\left(h_{\beta\rho}h_{\delta\alpha}h_{\omega\gamma}
           - h_{\beta\gamma}h_{\delta\rho}h_{\omega\alpha}\right) 
\end{array}
\end{equation}
\end{prop}

\bfseries\textit{Proof. }\normalfont
We want to find an expression for $h_{\alpha\beta;\gamma\delta}$. We start our calculations,

\[
\begin{array}{rcl}
h_{\alpha\beta;\gamma\delta}&=& h_{\alpha\gamma;\beta\delta}
+R_{\beta\gamma\alpha n;\delta}\\
&=& h_{\alpha\gamma;\delta\beta}+
g^{\rho\omega}\left(\overline{R}_{\beta\delta\gamma\rho}h_{\omega\alpha}
+\overline{R}_{\beta\delta\alpha\rho}h_{\omega\gamma}\right)
+R_{\beta\gamma\alpha n; \delta}\\
&=&\left(h_{\delta\gamma;\alpha}+R_{\alpha\delta\gamma n}\right)_{;\beta}\\
&& + g^{\rho\omega}\left(\overline{R}_{\beta\delta\gamma\rho}h_{\omega\alpha}
+\overline{R}_{\beta\delta\alpha\rho}h_{\omega\gamma}\right)
+R_{\beta\gamma\alpha n; \delta}\\
&=&h_{\delta\gamma;\alpha\beta}+R_{\alpha\delta\gamma n;\beta}\\
&& + g^{\rho\omega}\left(\overline{R}_{\beta\delta\gamma\rho}h_{\omega\alpha}
+\overline{R}_{\beta\delta\alpha\rho}h_{\omega\gamma}\right)
+R_{\beta\gamma\alpha n; \delta} ,
\end{array}
\]
then by Gauss equation,
\[
\begin{array}{rcl}
\overline{R}_{\beta\delta\gamma\rho}&=& 
R_{\beta\delta\gamma\rho}+h_{\beta\rho}h_{\delta\gamma}-h_{\beta\gamma}h_{\delta\rho}\\
\overline{R}_{\beta\delta\alpha\rho}&=& 
R_{\beta\delta\alpha\rho}+h_{\beta\rho}h_{\delta\alpha}-h_{\beta\alpha}h_{\delta\rho}
\end{array}
\]
from which we obtain
\begin{equation}
\label{Identity1}
\begin{array}{rcl}
h_{\alpha\beta;\gamma\delta}&=&
h_{\delta\gamma;\alpha\beta}+R_{\alpha\delta\gamma n;\beta}
+R_{\beta\gamma\alpha n; \delta}\\
&&+g^{\rho\omega}\left(R_{\beta\delta\gamma\rho}h_{\omega\alpha}+
R_{\beta\delta\alpha\rho}h_{\omega\gamma}\right)\\
&& g^{\rho\omega}(h_{\beta\rho}h_{\delta\gamma}h_{\omega\alpha}+
h_{\beta\rho}h_{\delta\alpha}h_{\omega\gamma}\\
&&\,\,\,\,\,\,\,\,\,\,\,\,-h_{\beta\gamma}h_{\delta\rho}h_{\omega\alpha}
-h_{\beta\alpha}h_{\delta\rho}h_{\omega\gamma}) .
\end{array}
\end{equation}

Using the following fact,
\begin{equation}
R_{\alpha\delta\gamma n;\beta}
= R_{\alpha\delta\gamma n|\beta}-h_{\alpha\beta}R_{n\delta\gamma n}
-h_{\delta\beta}R_{\alpha n\gamma n}
+g^{\rho\omega}h_{\beta\rho}R_{\alpha\delta\gamma\omega}
\end{equation}
and (\ref{Identity1}) we get
\[
\begin{array}{rcl}
h_{\alpha\beta;\gamma\delta}&=&
h_{\delta\gamma;\alpha\beta}+R_{\alpha\delta\gamma n|\beta}
+R_{\beta\gamma\alpha n|\delta}\\
&&-h_{\alpha\beta}R_{n \delta\gamma n}-h_{\delta\beta}R_{\alpha n\gamma n}
+g^{\rho\omega}h_{\beta\rho}R_{\alpha\delta\gamma\omega}\\
&&-h_{\beta\delta}R_{n \gamma\alpha n}-h_{\gamma\delta}R_{\beta n \alpha n}
+g^{\rho\omega}h_{\delta\rho}R_{\beta\gamma\alpha\omega}\\
&&+g^{\rho\omega}\left(R_{\beta\delta\gamma\rho}h_{\omega\alpha}+
R_{\beta\delta\alpha\rho}h_{\omega\gamma}\right)\\
&& g^{\rho\omega}(h_{\beta\rho}h_{\delta\gamma}h_{\omega\alpha}+
h_{\beta\rho}h_{\delta\alpha}h_{\omega\gamma}\\
&&\,\,\,\,\,\,\,\,\,\,\,\,-h_{\beta\gamma}h_{\delta\rho}h_{\omega\alpha}
-h_{\beta\alpha}h_{\delta\rho}h_{\omega\gamma}) .
\end{array}
\]

Contract in the indices $\gamma$ and $\delta$. Denoting by $H=g^{\alpha\beta}h_{\alpha\beta}$, 
by $\left|A\right|^2$ the norm of the tensor $h_{\alpha\beta}$, and using the 
Second Bianchi identity, produce
\begin{equation}
\begin{array}{rcl}
\Delta h_{\alpha\beta}&=& 
H_{\alpha\beta}+R_{\alpha n|\beta}+R_{\beta n|\alpha}-R_{\beta\alpha |n}+R_{n\alpha\beta n|n}\\
&&-h_{\alpha\beta}R_{nn}+HR_{n\alpha\beta n}-\\
&&+g^{\gamma\delta}g^{\rho\omega}h_{\beta\rho}R_{\alpha\delta\gamma\omega}
  +g^{\gamma\delta}g^{\rho\omega}h_{\delta\rho}R_{\beta\gamma\alpha\omega}\\
&&+Hg^{\rho\omega}h_{\beta\rho}h_{\omega\alpha}-|A|^2h_{\alpha\beta}\\
&&+g^{\gamma\delta}g^{\rho\omega}\left(R_{\beta\delta\gamma\rho}h_{\omega\alpha}
  +R_{\beta\delta\alpha\rho}h_{\omega\gamma}\right)\\
&&+g^{\gamma\delta}g^{\rho\omega}\left(h_{\beta\rho}h_{\delta\alpha}h_{\omega\gamma}
           - h_{\beta\gamma}h_{\delta\rho}h_{\omega\alpha}\right) .
\end{array}
\end{equation}

Finally, we use the fact
\begin{equation}
\frac{\partial}{\partial t}h_{\alpha\beta}=-R_{nn}h_{\alpha\beta}
+R_{\alpha n|\beta}+R_{\beta n|\alpha}-R_{\beta\alpha |n} ,
\end{equation}
to finally show the formula.

\hfill $\Box$

Contracting the previous formula in the indices $\alpha$ and $\beta$ yields,

\begin{cor}
\label{normalderivative1}
\begin{equation}
R_{\nu\nu |\nu}=-g^{\alpha \beta}\frac{\partial}{\partial t}h_{\alpha\beta}-HR_{\nu\nu}
=2\kappa\left[g^{\alpha\beta}R_{\alpha\beta}-R_{\nu\nu}\right].
\end{equation}
\end{cor}
\hfill $\Box$

\begin{prop}
\label{normalderivative2} 
\begin{equation}
\nabla_{\nu}R_{\alpha\beta}=R_{\nu\nu}h_{\alpha\beta} .
\end{equation}
\end{prop}

\noindent
\bfseries \textit{Proof. }\normalfont Let $\frac{\partial}{\partial \eta}$ be the outward normal
vector to $\partial M$ at $t=0$. Then, as a consequence of the Codazzi equations,
 $\frac{\partial}{\partial \eta}$ remains
normal to $\partial M$ throughout the flow. In this case we have that
\[
h_{\alpha\beta}=\frac{1}{2\left(g_{\eta\eta}\right)^{\frac{1}{2}}}\frac{\partial g_{\alpha\beta}}{\partial \eta} .
\]

A straightforward computation gives, 
\[
h_{\alpha\beta}'=R_{\nu\nu}h_{\alpha\beta}
-\frac{\partial R_{\alpha\beta}}{\partial \nu} .
\]

Recalling that $h_{\alpha\beta}=\kappa g_{\alpha\beta}$, we obtain
\[
-2\kappa R_{\alpha\beta}=R_{\nu\nu}h_{\alpha\beta}-\frac{\partial R_{\alpha\beta}}{\partial \nu} ,
\]
which can be written as,
\begin{equation}
\frac{\partial R_{\alpha\beta}}{\partial \nu}= 2\kappa R_{\alpha\beta}+R_{\nu\nu}h_{\alpha\beta} .
\end{equation}

On the other hand we have,
\[
\nabla_{\nu}R_{\alpha\beta}= \frac{\partial R_{\alpha\beta}}{\partial \nu}
-g^{\rho\omega}\left(R_{\alpha\rho}h_{\omega\beta}+R_{\beta\rho}h_{\omega\alpha}\right) ,
\]
which in the special case we are studying becomes,
\[
\begin{array}{rcl}
\nabla_{\nu}R_{\alpha\beta}&=& \frac{\partial R_{\alpha\beta}}{\partial \nu}
-\kappa g^{\rho\omega}\left(R_{\alpha\rho}g_{\omega\beta}+R_{\beta\rho}g_{\omega\alpha}\right)\\
&=& R_{\nu\nu}h_{\alpha\beta}+2\kappa R_{\alpha\beta}-2\kappa R_{\alpha\beta}\\
&=& R_{\nu\nu}h_{\alpha\beta} .
\end{array}
\]

\hfill $\Box$

\begin{rem}
The previous computation is valid even when the metric is not rotationally symmetric.
\end{rem}

From Corollary \ref{normalderivative1} and Proposition \ref{normalderivative2}
we get the following formula for the normal derivative of the scalar curvature

\begin{cor}
\begin{equation}
\frac{\partial}{\partial \nu}R=2\kappa g^{\alpha\beta}R_{\alpha\beta} .
\end{equation}
\end{cor}

\section{Proof of the Main Theorem}

We show that the estimates used in \cite{Ham82} are still valid in the case we are studying in this paper.
\subsection{Pinching Estimates}\label{Pinching}

The first consequence of the previous computations is the following result that we state without proof.

\begin{lem}
\label{positivityispreserved}
If $Ric>0$ at time $t=0$, it remains so as long as the solution to the flow exists.
\end{lem}

Using the results of the previous section we can show the following pinching estimates.
We adopt the following notation: $\mu\leq\lambda$ are the eigenvalues of $R_{\alpha\beta}$
and $\nu=R_{\nu\nu}$ on $\partial M$.

\begin{lem}
\label{lemmapinching}
Let $\eta>0$ be such that $\lambda<\left(1+\eta\right) \mu$ throughout the flow, and
assume that at $t=0$ we have $R_{ij}>\epsilon R g_{ij}$ for $\epsilon<\frac{1}{2\left(2+\eta\right)}$. 
This condition is preserved under the flow.
\end{lem}
\bfseries \textit{Proof. }\normalfont The equation satisfied in the interior by the 
tensor
\[
T_{ij}=\frac{R_{ij}}{R}-\epsilon g_{ij}
\]
are computed in \cite{Ham82}. All we have to see is what happens 
with the normal covariant derivative of $T$ at a point in the 
boundary when it gets a null eigenvalue. By the decomposition
of the Ricci tensor at the boundary
(consequence of Codazzi equations), we have to consider two cases separately: when
the eigenvector is normal to the boundary, and when it is tangent.

For the first case we must have 
\[
R_{\alpha\beta}\geq \epsilon Rg_{ij}\quad \mbox{and} \quad R_{\nu\nu}=\epsilon R .
\]  

If we compute the relevant part (it is clear that if $v\perp \partial M$, then
 $\left(\nabla_{\nu}T_{\alpha\beta}\right)v^{\alpha}v^{\beta}=0$),
 
\[
\begin{array}{rcl}
\nabla_{\nu}\left(\frac{R_{\nu\nu}}{R}-\epsilon\right)&=&
2\kappa g^{\alpha\beta}\frac{R_{\alpha\beta}}{R}-2\kappa \frac{R_{\nu\nu}}{R}
-\frac{1}{R^2}\left(2\kappa g^{\alpha\beta}R_{\alpha\beta}\right)R_{\nu\nu}\\
&=&2\kappa g^{\alpha\beta}\frac{R_{\alpha\beta}}{R}-2\kappa g^{\alpha\beta}
\frac{R_{\alpha\beta}}{R}\frac{R_{\nu\nu}}{R}-2\kappa\frac{R_{\nu\nu}}{R}\\
&\geq& 2\kappa g^{\alpha\beta}\frac{R_{\alpha\beta}}{R}-
 2\kappa g^{\alpha\beta}\frac{R_{\alpha\beta}}{R}\epsilon - 2\kappa\epsilon\\
&=& 2\kappa\left(1-\epsilon\right)g^{\alpha\beta}\frac{R_{\alpha\beta}}{R}-2\kappa\epsilon\\
&\geq& 2\kappa\left(1-\epsilon\right)g^{\alpha\beta}\left(\epsilon g_{\alpha\beta}\right)
-2\kappa\epsilon \\
&=& 4\kappa \left(1-\epsilon\right)\epsilon-2\kappa\epsilon\\
&=& 2\kappa\epsilon\left(1-2\epsilon\right)>0 \quad \mbox{if} \quad \epsilon\leq \frac{1}{2} .
\end{array}
\]

The second case is taken care of by the following computation (in this case
notice that if $v\in T_p\partial M$, then $\left(\nabla_{\nu}T_{\nu\nu}\right) v^n v^n= 0$
-because $v^n=0$). In this case we assume $R_{\alpha\beta}=\epsilon R g_{\alpha\beta}$
(here we use the hypothesis on the rotational symmetry of the metric)
and $R_{\nu\nu}\geq \epsilon R$.
\[
\begin{array}{rcl}
\nabla_{\nu}\left(\frac{R_{\alpha\beta}}{R}-\epsilon g_{\alpha\beta}\right)
&=& \frac{\nabla_{\nu}R_{\alpha\beta}}{R}
-\frac{1}{R^2}\frac{\partial R}{\partial \nu}R_{\alpha\beta}\\
&=& \frac{R{\nu\nu}}{R}h_{\alpha\beta}
-\frac{1}{R^2}2\kappa g^{\rho\sigma}R_{\rho\sigma}R_{\alpha\beta}\\
&=& \frac{R{\nu\nu}}{R}h_{\alpha\beta}
-2\kappa g^{\rho\sigma}\frac{R_{\rho\sigma}}{R}\frac{R_{\alpha\beta}}{R}\\
&\geq& \epsilon h_{\alpha\beta}
-2\kappa \left(2+\eta\right)\epsilon\cdot\epsilon g_{\alpha\beta}\\
&=& \kappa \epsilon g_{\alpha\beta}-4\kappa \epsilon^2 g_{\alpha\beta}
=\kappa\left[\epsilon-2\left(2+\eta\right)\epsilon^2\right]\\
&\geq& 0 \quad \mbox{if} \quad \epsilon\leq \frac{1}{2\left(2+\eta\right)} .
\end{array}
\]
\hfill $\Box$

As a corollary we get 
\begin{cor}
In the rotationally symmetric case if $R_{ij}\geq \epsilon Rg_{ij}$ at time $t=0$, with $\epsilon<\frac{1}{4}$,
then it remains so for all time.
\end{cor}
\bfseries\textit{Proof. }\normalfont Apply Lemma \ref{lemmapinching} with $\eta=0$.

\hfill $\Box$

\begin{cor}
The scalar curvature blows up in finite time.
\end{cor}
\bfseries \textit{Proof. }\normalfont  Notice that the scalar curvature satisfies the 
following differential inequality
\[\left\{
\begin{array}{l}
\frac{\partial R}{\partial t}\geq \Delta R + \epsilon^2 R^2\\
\frac{\partial R}{\partial \nu} \geq 0
\end{array}
\right.
\]
The Maximum Principle can be applied now to show the statement of the corollary.

\hfill $\Box$

\begin{rem}
From the proof of Lemma \ref{lemmapinching} we see that a sufficient condition to have the
pinching estimate is that the following holds. There is $\delta>0$ such that
\begin{equation}
\label{H-condition}
R_{\nu\nu}\geq \delta g^{\alpha\beta}R_{\alpha\beta} \quad\mbox{on}\quad \partial M\times \left(0,T\right)
\end{equation}

More exactly we have
\begin{prop}
\label{propoH-condition}
Assume inequality (\ref{H-condition}) holds as long as the solution to the 
Ricci flow exists. Then if at $t=0$ we have $R_{ij}\geq \epsilon R g_{ij}$ with $\epsilon\leq \frac{\delta}{2}$,
it remains so.
\end{prop}

{\bf\textit{Proof. }} The computations
in Lemma \ref{lemmapinching} the first (in time) null eigenvector of the
tensor
\[
T_{ij}=\frac{R_{ij}}{R}-\epsilon g_{ij}
\]
cannot happen in the normal direction. So if
$v$ is the first null eigenvector of $T_{ij}$ to occur, then $v^n=0$. We compute as follows,
\[
\begin{array}{rcl}
\nabla_{\nu}\left(\frac{R_{\alpha\beta}}{R}-\epsilon g_{\alpha\beta}\right)v^{\alpha}v^{\beta}
&\geq& \left[\delta g^{\rho\sigma}\frac{R{\rho\sigma}}{R}h_{\alpha\beta}
-2\kappa g^{\rho\sigma}\frac{R_{\rho\sigma}}{R}\frac{R_{\alpha\beta}}{R}\right]v^{\alpha}v^{\beta}\\
&\geq& \kappa g^{\rho\sigma}\frac{R_{\rho\sigma}}{R}\left[\delta g_{\alpha\beta}
-2\epsilon g_{\alpha\beta}\right]v^{\alpha}v^{\beta}\\
&\geq& 0 \quad \mbox{if} \quad \epsilon\leq \frac{\delta}{2} .
\end{array}
\]

\hfill $\Box$
\end{rem}
\subsection{Preserving pinching. }

We continue to use the following notation: $\lambda\geq \mu$ are the eigenvalues of $R_{\alpha\beta}$, and
$\nu=R_{\nu\nu}$. 
\begin{lem}
\label{unionofboundary} 
Let $\eta>0$ be such that for any $P\in \partial M$, we have
\[
\mbox{either}\quad \lambda\leq \left(1+\eta\right)\mu \quad
\mbox{or}\quad \lambda\leq \left(1+\eta\right)\nu .
\]

Then a pinching condition holds for $\epsilon\leq \frac{1}{2\left(2+\eta\right)}$.
\end{lem}

\bfseries\textit{Proof. }\normalfont
Choose $0<\epsilon < \frac{1}{2\left(2+\eta\right)}$ such that the pinching condition
$R_{ij}>\epsilon Rg_{ij}$ holds at $t=0$. Then if there is a time $t_0$ where this
pinching ceases to hold, there should be a point $p\in \partial M$ where the tensor
$T_{ij}=\frac{R_{ij}}{R}-\epsilon g_{ij}$ achieves its first zero eigenvalue. 
By the computations in Lemma \ref{lemmapinching} we know that the eigenvector corresponding
to this eigenvalue is tangent to the boundary. If at this point $\lambda<\left(1+\eta\right)\mu$
holds, the calculations in Lemma \ref{lemmapinching} show that  
\[
\nabla_{\nu}\left(\frac{R_{\alpha\beta}}{R}-\epsilon g_{\alpha\beta}\right)\geq 
\epsilon\kappa-2\kappa\left(2+\eta\right)\epsilon^2\geq 0.
\]

In the case that $\lambda < \left(1+\eta\right)\nu$ holds then $R_{\nu\nu}>2\epsilon g^{\alpha\beta}R_{\alpha\beta}$
holds, and at this point there cannot be a 0 eigenvalue of $T_{ij}$ by Proposition \ref{propoH-condition}. 

\hfill $\Box$

Consider the function 
\[
f=\frac{S}{R^2}=\frac{\lambda^2+\mu^2+\nu^2}{\left(\lambda+\mu+\nu\right)^2}
\]
then we have,

\begin{lem}
If there is $\delta>0$ such that $f\leq 1-\delta$, then there is $\eta>0$ such that the hypothesis of
Lemma \ref{unionofboundary} holds.
\end{lem}

\bfseries\textit{Proof. }\normalfont If the hypothesis of this Lemma holds, then we must have
\[
c\leq \frac{\mu}{\lambda}+\frac{\nu}{\lambda} < C,
\]
and the conclusion of the Lemma follows.

\hfill $\Box$

The previous two Lemmas show that 
if the pinching ceases to hold, it is because $f$ approaches 1 as $t\rightarrow T$.
Let us compute $\nabla_{\nu} f$ at $\partial M$.

\begin{lem}
We have
\[
\nabla_{\nu} f=2\frac{\kappa}{R^3}\left\{\nu\left(\lambda+\mu+\nu\right)\left[3\left(\lambda+\mu\right)-2\nu\right]
-2\left(\lambda+\mu\right)\left(\lambda^2+\mu^2+\nu^2\right)\right\}
\]
\end{lem}

A direct analysis of the previous expression shows the following,
\begin{lem}
\label{negativenormal}
There exists $\rho>0$ small enough such that if $\frac{\mu}{\lambda}+\frac{\nu}{\lambda}<\rho$ then
$\nabla_{\nu} f\leq 0$. Also if $\frac{\lambda}{\nu}\rightarrow 0$, the same conclusion holds.
\end{lem}

From \cite{Ham82} we borrow the following Lemma,

\begin{lem}
\label{evolutionofratio}
$f$ satisfies the following differential inequality
\[
\frac{\partial}{\partial t}f\leq \Delta f + u^k\partial_k f 
\]
where $u_k=\frac{2}{R}g^{kl}\partial_l R$.
\end{lem}

We show now that the maximum of $f$ cannot approach 1 as $t\rightarrow T$.
Indeed, if it does, we must have that either for every $\rho>0$ there is a time
$t$ where 
\[
\frac{\mu}{\lambda}+\frac{\nu}{\lambda}<\rho \quad \mbox{or}\quad \frac{\nu}{\lambda}\rightarrow \infty 
\]

But at such point $f$ cannot achieve a maximum, since this point must be in  $\partial M$ by Lemma
\ref{evolutionofratio}, and
Lemma \ref{negativenormal} would lead to a contradiction. Hence we have shown

\begin{thm}
There is an $\epsilon>0$ small enough, such that if $R_{ij}>\epsilon Rg_{ij}$ holds at time $t=0$, then it continues
to hold for all time.
\end{thm}

\hfill $\Box$


\subsection{Pinching the eigenvalues}

In this section we show that as it is the case in closed manifolds of positive Ricci curvature,
the eigenvalues of the Ricci tensor approach each other as the scalar curvature blows up.

\begin{thm}
\label{deltapinching}
We can find a $\delta>0$ and a constant $C$ depending only on the initial metric such that on $0\leq t<T$
we have 
\[
S-\frac{1}{3}R^2\leq CR^{2-\delta}.
\]
\end{thm}

\bfseries\textit{Proof. }\normalfont As we now have proved the pinching estimates, Hamilton's computations
for manifolds without boundary, carries over to the interior of oor manifold. If we define
\[
f=\frac{S}{R^{\gamma}}-\frac{1}{3}R^{\gamma-2},\quad \gamma=2-\delta,\quad \delta\leq 2\epsilon^2
\]
we know that $f$ satisfies an inequality
\[
\frac{\partial}{\partial t}f\leq \Delta f+u_k\partial_k f. 
\]

All we have to show then, to finally prove the Theorem is that for $\delta>0$ small enough, $\nabla_{\nu} f\leq 0$
outside a compact set of values of $\left(\frac{\lambda}{\nu},\frac{\mu}{\nu},1\right)$.
We do this in a series of Lemmas.

\hfill $\Box$

\begin{lem}
\[
\begin{array}{rcl}
\nabla_{\nu} f &=&
2\kappa\left\{\frac{\nu\left[3\left(\lambda+\mu\right)-2\nu\right]}{\left(\lambda+\mu+\nu\right)^{\gamma}}
-\gamma\frac{\left(\lambda^2+\mu^2+\nu^2\right)\left(\lambda+\mu\right)}{\left(\lambda+\mu+\nu\right)^{\gamma+1}}\right.\\
&& \left.+\delta\frac{1}{3}\left(\lambda+\mu+\nu\right)^{\gamma-3}\left(\lambda+\mu\right)\right\}
\end{array}
\]
\end{lem}

\begin{lem}
For $\delta>0$ small enough and $M>0$ big enough, if $\frac{\lambda}{\nu}+\frac{\mu}{\nu}\geq M$ then $\nabla_{\nu}f<0$.
\end{lem}

\bfseries\textit{Proof. }\normalfont
Assume for a moment that $\delta = 0$. Then the expression that determines the sign of $\nabla_{\nu} f$ is
\[
\left(\lambda+\mu+\nu\right)\left[3\left(\lambda+\mu\right)-2\nu\right]-2\left(\lambda^2+\mu^2+\nu^2\right)\left(\lambda+\mu\right).
\]

Dividing by $\nu$ and writing $x=\frac{\lambda}{\nu}, y=\frac{\lambda}{\mu}$, we must analyze the behavior of the 
\[
h\left(x,y\right)=\left(x+y+1\right)\left[3\left(x+y\right)-2\right]-2\left(x^2+y^2+1\right)\left(x+y\right).
\]

Notice that if $x=0$ or $y=0$, there is $\rho>0$ such that $h\left(x,y\right)\leq -\rho$. Now we show that for
a rectangle $Q$ big enough with sides parallel to the coordinate axes, we must have that if
$\left(x,y\right)\in \mathbf{R}^2\setminus Q$, then $h\left(x,y\right)<0$. 
Indeed, if we maximize the function $h$ subject to the constraint $x+y=\epsilon$, then we obtain that the
maximum is reached when $x=\frac{\epsilon}{2}=y$. Evaluating at this point we obtain
\[
h\left(\frac{\epsilon}{2},\frac{\epsilon}{2}\right)=
-\epsilon^3+3\epsilon^2-\epsilon-2
\]
it is easy to see that if $\epsilon>3$ and $\epsilon<1$, then $h\left(\frac{\epsilon}{2},\frac{\epsilon}{2}\right)<0$. This
together with the previous observation shows the assertion. 

Observe that the term
\[
\frac{\left(\lambda+\mu+\nu\right)^{\gamma}\left(\lambda+\mu\right)}{\left(\lambda+\mu+\nu\right)^3}
\]
is bounded if $R\geq \rho>0$. This show that if $\delta>0$ is small enough, then $\nabla_{\nu}f<0$.
This also finishes the proof of Theorem \ref{deltapinching}.

\hfill $\Box$
\subsection{The gradient of the Scalar Curvature. }

Now we can use a compactness argument as in \cite{Ham95} to show the following
\begin{thm}
For every $\theta>0$ we can find s constant $C\left(\theta\right)$ depending only on $\theta$ and the 
initial value of the metric, such that on $0\leq t<T$ we have
\[
\max_{t\leq \tau}\max_P\left|D Rm\left(P,t\right)\right|
\leq \theta \max_{t\leq\tau}\max_P\left|Rm\left(P,\tau\right)\right|^3+C\left(\theta\right) 
\]
\end{thm}

Of course, to produce this
compactness argument we need
to bound derivatives of the curvature in terms of the curvature.
The rest of the proof of long time existence and exponential convergence then follows the same arguments as in
\cite{Ham82}, again as long as we learn how to bound derivatives of the curvature in terms of bounds in
the curvature, and this is the purpose of the last section of this work.

\section{Bounding derivatives of the the curvature}\label{Derivative}

In this section we sketch a 
procedure to produce bounds on the derivatives of the curvature tensor from bounds on the curvature. We skip
the first order derivatives since the computations of Theorem 7.1 in \cite{Ham95}
(see also Theorem 6.1 in \cite{Cortissoz}), can be easily adapted to
our case and show
how to produce bounds on certain second derivatives, hoping that
the reader will find easy to convince himself that this method 
extends to bound any number of derivatives. 

Fist of all, fix a collar of the boundary  and
fix Fermi coordinates $\left(x^1,\dots, x^n\right)$ with respect to the initial metric, where 
$x^n$ is the distance to the boundary. The vector fields $\partial_{\alpha}$, $\alpha=1,\dots, n-1$
remain tangent to the boundary (of course when restricted to the boundary) whereas the vector field $\partial_n$ remains normal (recall
that we are working with a weakly umbilic boundary it remains).

By the interior derivative estimates of Shi, on the inner boundary of the collar for $t>0$ we 
can assume bounds on all the derivatives of the curvature in terms of a bound of the curvature. 
All these said, we set ourselves to the task of estimating derivatives of the form $\nabla_i\nabla_{\nu} Ric$, $i=1,\dots, n$
where we denote by $\nu:=\frac{1}{g_{nn}}\partial_n$, which 
coincides with the outward unit normal when restricted to the boundary.
From now on we assume bounds $\left|Ric\right|\leq M_0$ and
$\left|\nabla Ric\right|\leq M_1$ on $\overline{M}\times \left[0,T\right]$, $T<\infty$.

The evolution equations for the Ricci tensor are given by
\begin{equation}
\label{important1}
\left\{
\begin{array}{l}
\frac{\partial}{\partial t}R_{\alpha\beta}=\Delta R_{\alpha\beta}-Q_{\alpha\beta}
\quad \mbox{in} \quad M\times\left(0,T\right)\\
\nabla_{\nu}R_{\alpha\beta}=\kappa R_{\nu\nu}g_{\alpha\beta} \quad \mbox{on} \quad \partial M \times \left(0,T\right)
\end{array}
\right.
\end{equation}
\begin{equation}
\label{important2}
\left\{
\begin{array}{l}
\frac{\partial}{\partial t}R_{\nu\nu}=\Delta R_{\nu\nu}-Q_{\nu\nu}+2\left(R_{\nu\nu}\right)^2
\quad \mbox{in} \quad M\times\left(0,T\right)\\
\nabla_{\nu}R_{\nu\nu}=2\kappa\left[g^{\alpha\beta}R_{\alpha\beta}-R_{\nu\nu}\right] \quad
\mbox{on} \quad \partial M\times\left(0,T\right)
\end{array}
\right.
\end{equation}

Borrowing from Lemma 13.1 in \cite{Ham82} we find that the first covariant derivative of the
Ricci tensor satisfies the following system of equations
\begin{prop}
\begin{equation}
\label{System}
\left\{
\begin{array}{l}
\frac{\partial}{\partial t}\left(\nabla_{\nu}Ric\right)=
\Delta \left(\nabla_{\nu} Ric\right) +Ric*\nabla_{\nu} Ric \quad \mbox{in}\quad M\times\left(0,T\right)\\
\nabla_{\nu} Ric = g*Ric \quad\mbox{on}\quad \partial M\times\left(0,T\right)
\end{array}
\right.
\end{equation}
\end{prop}
as in \cite{Ham82}, if $A$ and $B$ are two tensors we write $A*B$ for any linear combination
of tensors formed by contraction on $A_{i\dots j}B_{k\dots l}$ using the metric $g^{ik}$.
Notice that in the expession $g*Ric$ we have absorbed the constant $\kappa$.

We have to compute the Laplacian when acting on 3-tensors. We have the following formula
\begin{prop}
\label{laplacian3tensors}
\begin{equation}
\begin{array}{rcl}
\nabla_j\nabla_i \Psi_{lmk}&=&\frac{\partial^2}{\partial x_i\partial x_j}\Psi_{lmk}\\
&&-\Psi_{pmk}\frac{\partial}{\partial x_j}\left(\Gamma_{il}^p\right)
-\Psi_{lpk}\frac{\partial}{\partial x_j}\left(\Gamma_{im}^p\right)
-\Psi_{lmp}\frac{\partial}{\partial x_j}\left(\Gamma_{ik}^p\right)\\
&&-\Gamma_{il}^p\Psi_{pmk,j}-\Gamma_{im}^p\Psi_{lpk,j}-\Gamma_{ik}^p\Psi_{lmp,j} 
\end{array}
\end{equation}
\end{prop}

Therefore, taking $\Psi=Ric$ in Proposition \ref{laplacian3tensors}, system (\ref{System}), becomes
\[
\left\{
\begin{array}{l}
\frac{\partial}{\partial t}\nabla_{\nu}Ric -\Delta_e \nabla_{\nu}Ric = \mathcal{W} \quad \mbox{in} \quad M\times \left(0,T\right)\\
\nabla_{\nu} Ric = g*Ric \quad \mbox{on} \quad \partial M\times \left(0,T\right)
\end{array}
\right.
\]
where 
\[
\Delta_e = g^{ij}\frac{\partial^2}{\partial x^i\partial x^j} \quad \mbox{and} 
\]
\[
\mathcal{W}=\Gamma*\nabla^2 Ric+Ric*\nabla Ric+\partial \Gamma*\nabla Ric.
\]

Let $F_{ij}$ be the extension of the boundary quantity that appears in 
(\ref{important1}) and (\ref{important2}). Then we have the following formulae,
\[
\left\{
\begin{array}{c}
\frac{\partial F}{\partial t}= Ric * Ric + g*\nabla^2 Ric\\
\Delta F= g*\nabla^2 Ric
\end{array}
\right.
\]

Hence, if we define $U=\nabla Ric - F$, we find that it satisfies the following system
\[
\left\{
\begin{array}{l}
\frac{\partial U}{\partial t}-\Delta_e U=\overline{\mathcal{W}} \quad\mbox{in}\quad M\times\left(0,T\right)\\
U=0 \quad \mbox{on}\quad \partial M \times \left(0,T\right)
\end{array}
\right.
\]
where $\overline{\mathcal{W}}=\mathcal{W}-\frac{\partial F}{\partial t}-g*\nabla^2 Ric$. 

Via the Ricci flow, and taking into account that we are considering a bounded time interval,
it can be shown that $\left|\partial^i \Gamma\right|\leq C\left|\nabla^{i+1} Ric\right|$, 
and hence we can bound 
 $\overline{\mathcal{W}}$ in terms of bounds on $\left|\nabla^i Ric\right|$ ($i=1,2$).
 
By the local interior estimates of Shi, we can assume that $U=0$ on the part of the
boundary of the collar which is contained in the interior of the manifold. We must
point out that by assuming $U=0$ we are introducing in the estimates for the bounds
terms depending on the size of the collar of the boundary (at time $t=0$) where we are estimating.

By the theory of the
first boundary value problem, we can represent $U$ as
\begin{equation}
\label{U}
U=\int_0^\theta\left(\int\Gamma \overline{\mathcal{W}}\right)\,d\tau +\int \Gamma \left(\nabla Ric+g*Ric\right)\left.\right|_{t=0}
\end{equation} 
where $\Gamma$ is the fundamental solution of $\Delta_e$. The integral sign without specified
region of integration refers to spatial integration.

Notice that $\nabla U = \nabla^2 Ric$. Therefore, if $M_2=\max \left|\nabla^2 Ric\right|$ we get
from (\ref{U})
\[
M_2\leq M_2\int_0^\theta \left(\int \nabla\Gamma\right)+\frac{M_1\cdot M_0}{\sqrt{\theta}}
\]
here $M_i=\max\left|\nabla^i Ric\right|$, $i=0,1$

By making $\theta>0$ small we obtain
\[
M_2 \leq \frac{M_1M_0}{\sqrt{\theta}}
\]

The smallness of $\theta$ depends on Lipszchitz bounds on the coefficients of $\Delta_e$ which
in turn only depend on $M_0$ and $M_1$.


\end{document}